\documentclass[11pt]{article}
\usepackage{latexsym}
\usepackage{amsfonts}

\newcommand{\lf}{\lfloor}
\newcommand{\rf}{\rfloor}

\newcommand{\hsp}{\hspace{\parindent}}

\newcommand{\de}{\delta}
\newcommand{\eps}{\epsilon}
\newcommand{\la}{\lambda}
\newcommand{\sig}{\sigma}

\newcommand{\La}{\Lambda}
\newcommand{\Sig}{\Sigma}

\newcommand{\NN}{{\Bbb N}}
\newcommand{\RR}{{\Bbb R}}
\newcommand{\QQ}{{\Bbb Q}}
\newcommand{\ZZ}{{\Bbb Z}}

\newcommand{\sA}{{\cal A}}

\newcommand{\sS}{{\cal S}}

\newcommand{\sP}{{\cal P}}

\newcommand{\bc}{{\bf c}}
\newcommand{\bd}{{\bf d}}

\newcommand{\bp}{{\bf p}}

\newcommand{\bx}{{\bf x}}
\newcommand{\by}{{\bf y}}

\newcommand{\beql}[1]{\begin{equation}\label{#1}}
\newcommand{\eeq}{\end{equation}}

\catcode`\@=11
\renewcommand{\section}{
        \setcounter{equation}{0}
        \@startsection {section}{1}{\z@}{-3.5ex plus -1ex minus
        -.2ex}{2.3ex plus .2ex}{\large\bf}
        }
\catcode`@=12
\makeatletter
\def\eqalignno#1{\displ@y \ta {\bf s} kip\@centering
  \halign to\displaywidth{\hfil$\@lign\displaystyle{##}$\ta {\bf s} kip\z@skip
    & $\@lign\displaystyle{{}##}$\hfil\ta {\bf s} kip\@centering
    & \llap{$\@lign##$}\ta {\bf s} kip\z@skip\crcr
    #1\crcr}}
\makeatother

\makeatletter
\def\@sect#1#2#3#4#5#6[#7]#8{\ifnum #2>\c@secnumdepth
     \def\@svsec{}\else 
     \refstepcounter{#1}\edef\@svsec{\csname the#1\endcsname.\hskip .75em }\fi
     \@tempskipa #5\relax
      \ifdim \@tempskipa>\z@ 
        \begingroup #6\relax
          \@hangfrom{\hskip #3\relax\@svsec}{\interlinepenalty \@M #8\par}%
        \endgroup
       \csname #1mark\endcsname{#7}\addcontentsline
         {toc}{#1}{\ifnum #2>\c@secnumdepth \else
                      \protect\numberline{\csname the#1\endcsname}\fi
                    #7}\else
        \def\@svsechd{#6\hskip #3\@svsec #8\csname #1mark\endcsname
                      {#7}\addcontentsline
                           {toc}{#1}{\ifnum #2>\c@secnumdepth \else
                             \protect\numberline{\csname the#1\endcsname}\fi
                       #7}}\fi
     \@xsect{#5}}
\def\@theorem#1#2{\it \trivlist \item[\hskip \labelsep{\bf #1\ #2.}]}

\makeatother
\begin{document}

\newtheorem{Theorem}{Theorem}[section]
\newtheorem{Proposition}{Proposition}[section]
\newtheorem{Corollary}{Corollary}[section]
\newtheorem{Lemma}{Lemma}[section]
\newtheorem{defi}{Definition}[section]

\begin{center}
{\large {\bf Local Complexity of Delone Sets and Crystallinity}} \bigskip \\
{\large {\em Jeffrey C. Lagarias}} \\ \smallskip
AT\&T Labs -- Research \\
Florham Park, New Jersey 07932 \\
\vspace*{1\baselineskip}
{\large {\em Peter A. B. Pleasants}} \\ \smallskip
Dept.\ of Mathematics and Computing Science \\
University of the South Pacific \\
Suva, Fiji \\
\vspace{2\baselineskip}
(May 5, 2001 version) \\
\vspace{2\baselineskip}
{\em Abstract}
\end{center}

This paper characterizes when a Delone set $X\/$ in
$\RR^n$ is an ideal crystal in terms of restrictions
on the number of its local patches of a given size
or on the hetereogeneity of their distribution.
For a Delone set $X\/$, let $N_X (T)$ count the number
of translation-inequivalent patches of radius $T\/$
in $X\/$ and let $M_X(T)$ be the minimum radius such that
every closed ball of radius $M_X(T)$ contains the center
of a patch of every one of these kinds.  
We show that for each of these functions there is a
``gap in the spectrum'' of possible growth rates between
being bounded and having linear growth, and that having
linear growth is equivalent to $X\/$ being an ideal crystal.

Explicitly, for $N_X(T)$,  if $R\/$ is the covering
radius of $X\/$ then
either $N_X(T)$ is bounded or $N_X (T) \ge T/2R$ for all $T>0$.
The constant $1/2R\/$ in this bound is best possible in all dimensions.

For $M_X(T)$, either $M_X(T)$ is bounded or $M_X (T) \ge T/3$ for all $T>0$.
Examples show  that the constant $1/3$ in this bound cannot be replaced 
by any number exceeding $1/2$.
We also show that every aperiodic Delone set $X\/$ has $M_X(T)\ge c(n)\,T\/$
for all $T>0$, for a certain constant
 $c(n)$ which depends on the dimension $n\/$ of $X\/$ and
is $>1/3$ when $n>1$.

\vspace*{1.5\baselineskip}


\newpage

\begin{center}
{\large {\bf Local Complexity of Delone Sets and Crystallinity }} \bigskip \\
{\large {\em Jeffrey C. Lagarias}} \\ \smallskip
AT\&T Labs -- Research \\
Florham Park, New Jersey 07932 \\
jcl@research.att.com \bigskip \\
\vspace*{1\baselineskip}
{\large {\em Peter A. B. Pleasants}} \\ \smallskip
Dept.\ of Mathematics and Computing Science \\
University of the South Pacific \\
Suva, Fiji \\
pleasants\_p@usp.ac.fj \\
\vspace{1\baselineskip}
\vspace{2\baselineskip}
\end{center}
\setlength{\baselineskip}{1.5\baselineskip}

%
%

\section{Introduction}
\hsp
The discovery in 1984 of quasicrstalline materials, which are strongly
ordered aperiodic structures, has generated renewed interest in precisely
delineating the boundary between discrete sets that are fully periodic
(ideal crystals) and those with less global order, see Radin~\cite{Ra91}
and several of the articles in \cite{Moo97b}, edited by Moody.
This paper studies such questions in discrete geometry,
investigating how strong the restrictions on local features of 
discrete point sets in $\RR^n$ must be to enforce crystallinity.
Such questions originally arose in geometric crystallography,
see \cite{En93}.
The questions we study also have connections
with the ergodic theory of $\RR^n$-actions \cite{LP99};
there are parallel questions in the combinatorics of words,
in one or several dimensions, which we briefly consider. 

We study the patch-counting function $N_X(T)$ and the
repetitivity function $M_X(T)$ of Delone sets $X$, where
patches are identified only up to translation equivalence.
(These concepts were studied in \cite{LP99} and are defined in \S2.)
We show that both functions have a ``gap in the spectrum'' of their
possible growth rates: either they are bounded, which happens precisely
when $X\/$ is an ideal crystal, or they grow at least linearly in
the radius $T\/$ of the patch.  Growth rates like $\log T$ or $\sqrt T$ 
are impossible for these functions.  These results can be viewed as 
characterizations of ideal crystals: when $X\/$ is an ideal crystal both
functions are bounded, but nevertheless sufficiently slow linear growth
of either of these functions is enough to characterize ideal crystals.
They can alternatively be viewed as giving linear lower
bounds on the growth rates of these functions for Delone sets $X$
that are not ideal crystals.

The problem of characterizing ideal crystals in terms of restrictions
on allowed types of local patches was first studied for ``regular point
systems'' in 1976 by Delone~et~al.\ \cite{DDSG76} for isometry-equivalence
classes of local patches.
They identified a radius $T\/$ (depending on the dimension and the Delone
constants of $X\/$) such that every Delone set $X\/$ whose patches of
radius $T\/$ are all isometry-equivalent is a regular point system
and therefore an ideal crystal.
Dolbilin~et~al.\ \cite{DLS98} gave an extension of this result to
``multiregular point systems'', which comprise all ideal crystals:
for each $k \ge 1$ they gave a radius $T$ (depending on the dimension,
the Delone constants of $X\/$ and on $k$) such that if there are at
most $k$ isometry-equivalence types of patches of this radius, then
$X$ is an ideal crystal with at most $k$ isometry-equivalence types
of patches of {\em any\/} radius.
This gives an effective way of deciding, for a given $k$, whether $X\/$
has degree of regularity $\le k$.
It does not, however, give a way of deciding whether a given set $X\/$ has
some finite degree of regularity; that is, whether $X\/$ is an ideal crystal.
This latter drawback is unavoidable, because with any finite sample of $X\/$
there is no way to tell that it is not part of an ideal crystal whose unit
cell is larger than the size of the sample.
However, the dependence of $T\/$ on $k\/$ in \cite{DLS98} is linear, and it
follows that there is a constant $c\/$ (depending only on the dimension and
the Delone constants) such that if there is {\em some} $T\/$ with at most
$cT\/$ isometry patches of radius $T\/$ then $X\/$ is an ideal crystal.

For the patch-counting function $N_X(T)$, our result (Theorem~\ref{th21})
is an analogue of this small linear-growth rate characterization of ideal
crystals, in which patches are identified only up to translation-equivalence.
This is a finer equivalence relation on patches than isometry-equivalence,
and stronger bounds hold.
The corresponding coefficient of $T\/$ we obtain for translation-equivalence
types of patches (which gives a growth rate of $N_X(T)$ in the gap between
crystals and non-crystals) depends only on the covering radius of $X$,
not on its packing radius or dimension, and we show that it is best
possible in all dimensions. 
We conjecture that certain super-linear polynomial growth rates of
$N_X(T)$, while not guaranteeing that $X\/$ is an ideal crystal,
are sufficient to ensure that $X\/$ has some periods---more precisely
that if, for a sufficiently small $c$, $N_X(T)<cT^{n-j+1}$ for all
large $T\/$ then $X\/$ has $j\/$ independent periods---but no cases
other than $j=n$ have yet been proved.

For the repetitivity function $M_X(T)$ (which is scale-invariant)
the coefficient of $T\/$ we obtain in Theorem~\ref{th22} for a
linear growth rate in the gap between crystals and non-crystals
is the absolute constant $\frac{1}{3}$, independent of dimension.
We show by examples that this cannot be improved beyond $\frac{1}{2}$.
In Theorem~\ref{th23} we also show that a coefficient slightly larger
than $\frac{1}{3}$ (but depending on the dimension $n$) in the linear
growth rate of $M_X(T)$ is enough to ensure that $X\/$ has a non-zero
period when $n > 1$.
This coefficient depends on the packing-covering constant $\kappa(n)$
for $n$-dimensional Delone sets, a concept originally studied
by Ryshkov~\cite{R74}, which is defined in \S2.
In the final section we study various properties of this constant.

The results of this paper are analogous to results in the
combinatorics of words (symbolic dynamics) on the lattice $\ZZ^n$.
They can even be viewed as generalizing such results, because 
symbolic words on $\ZZ^n$ can be encoded as Delone sets by placing
points near the lattice points in $\ZZ^n$ using small dispacements
to distinguish the symbol types.
A well known result in the combinatorics of words on the lattice
$\ZZ$---the one-dimensional case---given in Morse and Hedlund
\cite{MH38} is that if $N_S(m)$ is the number of different words
of length $m\/$ in a two-sided infinite sequence $S\/$ in which
$A\/$ distinct symbols occur then either $S\/$ is periodic and
$N_S(m)\le P\/$ (the length of the period) for all $m\/$ or
\beql{Nforall}
N_S(m)\ge m+A-1\mbox{ for all $m$.}  
\eeq
Associated to any such symbol sequence is a dynamical system
with a $\ZZ$-action (the closure of its orbit under the shift)
and this result can be interpreted as a condition for this
dynamical system to be finite.
In the context of higher dimensional symbolic dynamical systems, with
symbols on the lattice $\ZZ^n$ and a corresponding $\ZZ^n$-action,
the question has also been raised to what extent growth restrictions
on the number of local symbol patterns (in rectangular patches)
enforce periodicity, see for example Sander and Tijdeman
\cite{ST99a,ST99b, ST00} and Berth\'e and Vuillon \cite{BV00,BV99}
for a discussion of the two dimensional case.
Symbolic dynamics analogues  of the patch-counting function are called
``permutation numbers'' or ``complexity functions''
in this context, cf.~\cite{Fe96,Lo83}.
The approach of Theorem~\ref{th21} can be used to show that
a low enough complexity bound on the number of rectangular
patches of such systems enforces full-dimensional periodicity.
For square (or, more generally, cubic) patches, however,
one can directly obtain an optimal result in $n\/$ dimensions,
given as Theorem~\ref{th41} in \S4; the proof is a straighforward
extension of the proof of Morse and Hedlund~\cite[Theorem~7.4]{MH38}
of (\ref{Nforall}) in one dimension.

The analogue in symbolic dynamics of the repetitivity function
is the {\em recurrence function} $R_S(m)$, introduced by
Morse and Hedlund~\cite{MH38,MH40} in 1938.
Given a two-sided infinite symbol sequence $S$, $R_S(m)$ is
the shortest length such that every word of length $R_S(m)$
contains a copy of every word of length $m$ that occurs in $S$. 
To aid comparison with our repetitivity function we shall describe
results in terms of a related function $M_S(m)=R_S(m)-m$, which
represents the maximum distance between the leading symbols of
any two successive identical words of length $m$.
Put in terms of $M_S(m)$, Morse and Hedlund showed
\cite[Theorem~7.5]{MH38} that for any aperiodic one-dimensional
repetitive sequence $S\/$ in which $A\/$ distinct symbols occur
\beql{Rforall}
M_S(m)\ge m+A-1\mbox{ for all $m\ge1$}
\eeq
and \cite[p.2]{MH40} that
\beql{limsup}
\limsup_{m \to \infty} \frac{M_S(m)}{m} \geq \tau + 1,
\eeq
where $\tau = \frac{1}{2}(1 + \sqrt{5})$ is the golden ratio.
Moreover (\ref{limsup}) holds with equality when $S\/$ is
the Fibonacci sequence, so represents an optimal
``gap in the spectrum'' result.
This suggests the possibility that the limiting constant $\tau+1$
in the analogous Delone set inequality on the right of (\ref{604})
in Theorem~\ref{th61} may also be optimal.
The inequality (\ref{Rforall}) follows from (\ref{Nforall}) and the
trivial estimate $M_S(m)\ge N_S(m)$ (the analogue for symbol sequences
of Theorem~4.1 of \cite{LP99}) yet the coefficient of $m\/$ in
(\ref{Rforall}) is 1, which is significantly larger than the $\frac{1}{3}$
in our Theorem~\ref{th22}, the analogous result for Delone sets.
This seems largely due to $m\/$ being restricted to integer values
in the symbolic case.

Some of the bounds in this paper depend on the 
$n$-dimensional Delone packing-covering constant $\kappa(n)$,
defined in \S2, which was introduced and studied by
Ryshkov~\cite{R74} in 1974.
In the final section we present
basic properties of this constant and the related 
$n$-dimensional lattice packing-covering
constant $\kappa_L(n)$, also introduced by Ryshkov. 
Most of the results in that section are not new, but we recall them
for completeness and 
in order to raise a question concerning the possible unboundedness
of $\kappa_L(n)$ as $n \to \infty$.
This question is of interest, because in any dimension $n\/$ with
$\kappa_L(n) > 2$ no lattice packing can be a densest sphere packing.

%
%

\section{Statements of results}\label{defs}
\hsp
\begin{defi}~\label{defs.1}
{\em 
A {\em Delone set}, or $(r,R)${\em-set}, is a discrete set $X\/$ in
$\RR^n$ that is {\em uniformly discrete\/} (i.e.\ its packing radius
$r\/$ by equal balls is positive) and {\em relatively dense\/}
(i.e.\ its covering radius $R\/$ by equal balls is finite). We
call the values $r$ and $R$ the {\em Delone constants} of the
set $X$.
}
\end{defi}

Clearly $r \le R$, with equality only when $X\/$ is a one-dimensional set
of equally spaced points.  Our definition of the Delone constant $r\/$
differs slightly from other authors (see \cite{DDSG76,DLS98,R74},
for example) who take for this value
 the infimum $2r\/$ of the inter-point distances
in place of the packing radius $r\/$.

\begin{defi}~\label{defs.2}
{\em 
A Delone set $X\/$ in $\RR^n$ is an  {\em ideal crystal\/} if it has
a full rank lattice of translation symmetries, i.e.\ $X=\La+F$, where
$\La$ is a full rank lattice in $\RR^n$ and $F\/$ is a finite set.
}
\end{defi}

Note that an {\em aperiodic set} is one with no global translation symmetries,
while a {\em non-crystalline set} may have some translation symmetries,
but not a full rank set of them.

\begin{defi}~\label{defs.3}
{\em
For a Delone set $X$, a {\em T-patch\/} centered at the 
point $\bx\in X$ is
$$
\sP_X (\bx ; T) : = X \cap {B}( \bx ; T)~,
$$
where ${B}( \bx ; T)$ is the open ball with center $\bx$ 
and radius $T$, and
the {\em patch-counting function} $N_X (T)$ is the number 
(possibly infinite)
of translation-inequivalent $T$-patches centered at points $\bx$ of $X$.
We use the notation $\sP(\bx;T)\sim\sP(\by;T)$ to mean translation equivalence
of patches, i.e.\ that $\sP(\by;T)=\sP(\bx;T)+\by-\bx$.
}
\end{defi}

The function $N_X(T)$ is a non-decreasing function of $T\/$ and $N_X(T)=1$ for
$T\le 2r$. It is
possible  that $N_X(T)$ may be infinite for large $T$.  

\begin{defi}~\label{defs.3a}
{\em
A Delone set $X\/$ with $N_X(T)$ finite for all $T\/$ is said 
to have {\em finite local complexity}.
}
\end{defi}

Delone sets with finite local complexity are studied in
\cite{La99a,La99b,LP99}, where they are called Delone sets of finite type.
The growth rate of the patch-counting function $N_X (T)$ provides a
quantitative  measure of the complexity of $X$.   The slowest possible
growth rate for $N_X(T)$ is to be eventually constant, which occurs
when $X\/$ is an ideal crystal, as we review in \S\ref{Nbdd}.
Our first main result is that except for ideal crystals the slowest
possible growth rate is at least linear in $T$.
Dolbilin~et~al.\ \cite[Theorem~1.3]{DLS98} have already shown, 
in the context of counting isometry classes of patches, 
that $N_X(T)<c(n,r,R)T\/$
for a single value of $T\/$ implies that $X\/$ is a crystal, where
$$
c(n, r, R) = \frac {1}{2( n^2 + 1) R \log_2 (\frac{R}{r} + 2)},
$$
and this holds {\em a fortiori\/} for our function $N_X(T)$ which
counts the more numerous translation classes of patches.  But with
our finer classification of patches we can obtain a larger value of the
constant, which is  independent of the dimension.

\begin{Theorem}~\label{th21}
If a Delone set $X$ in $\RR^n$ with covering radius
$R$ has a single value of $T > 0$ such that
\beql{Ncrystal}
N_X(T) < \frac{T}{2R}\,,
\eeq
then $X$ is an ideal crystal.
\end{Theorem}

The constant $\frac{1}{2R}$ here is optimal, in the sense that
for any $c > \frac{1}{2R}$ there are Delone sets $X\in\RR^n$
with  $N_X(T) < c T$, for some $T>0$, that are not ideal crystals,
as we show in \S\ref{examples}.

In \cite[Conjecture~2.2]{LP99} we put forward the following conjecture, 
which says
that slow growth of $N_X(T)$ implies that $X\/$ has many independent periods:

\paragraph{Period Conjecture.}{\em For each integer $j=1,\ldots,n$ there
is a positive constant $c_j(n,r,R)$ such that any Delone set $X$
in $\RR^n$ with Delone constants $r,R$ that satisfies
$$
N_X(T)<c_j(n,r,R)T^{n-j+1}\mbox{ for all $T>T_0(X)$}
$$
has $j\/$ linearly independent periods.}

Weaker forms of this conjecture would be to allow $c_j$ to depend
on $X\/$ or merely to assert that $T^{j-n-1}N_X(T)\to0$ as
$T\to\infty$ implies that $X\/$ has $j\/$ linearly independent periods.
Theorem~\ref{th21} is a strong form of the case $j=n\/$ of the
conjecture.  The only other case that has been proved to date is
a strong form of the case $n=2$, $j=1$ in the context of doubly
infinite arrays of symbols instead of Delone sets \cite{EMK}.

The basis of the conjecture is the feeling that a pattern with less
than $j\/$ independent periods should have ``$(j-1)$-dimensional
recognizable features'' at arbitrarily large scales, causing $N_X(T)$
to be of order at least $T^{n-j+1}$ for some arbitrarily large values
of $T$; but concrete evidence for it is limited.
We formulate it in order to focus attention on the problem,
rather than to strongly assert its truth.
The formulation above is the strongest consistent with
currently known constraints.
In particular, the conjecture cannot be strengthened to assert
the conclusion when the inequality is satisfied for only a single
value of $T$, as Theorem~\ref{th21} does: in \cite[Theorem~2.2]{LP99}
we give examples of aperiodic sets $X\subset\RR^n$, for any $n\ge3$
and $\eps>0$, for which there exist arbitrarily large radii $T\/$
with $N_X(T)<T^{\lceil(n+1)/2\rceil+\eps}$.  
Also there is no corresponding conjecture (except possibly for
$n=2$, $j=1$) generalizing the result of
\cite{DLS98} for the function that counts patches up to isometry: 
\cite{DP00} gives examples of aperiodic sets $X\subset\RR^n$, 
for $n\ge3$, whose
isometry patch-counting functions are $O(T^{1+\eps})$ for any $\eps>0$.

We next consider sets that have 
restrictions on the distribution of patches of a given type.

\begin{defi}~ \label{defs.4}
{\em
For a Delone set $X$, the {\em repetitivity function} $M_X (T)$ is the
least $M\/$ (possibly infinite) such that every closed ball $\bar B\/$
of radius $M\/$ contains the center of a $T$-patch of every kind that
occurs in $X$.  That is, for every $T\/$-patch $\sP$ of $X$, $\bar B\/$
contains a point of $X\/$ which is the center of a $T\/$-patch of $X\/$
that is a translate of $\sP$.  (This $T\/$-patch may extend outside the
ball $\bar B\/$.) 
} 
\end{defi}

Another way of expressing this definition is to say that $M_X(T)$ is
the largest of the covering radii of the sets of centers of patches
of $X\/$ translation-equivalent to $\sP$, taken over all $T$-patches
$\sP$ of $X\/$.  The function $M_X(T)$ is a non-decreasing function
of $T\/$ and $M_X(T)=R\/$ for $T\le 2r$.  For large $T$, $M_X(T)$ may
become infinite.  

\begin{defi}~\label{defs.4a}
{\em 
(i) A Delone set $X\/$
with $M_X(T)$ finite for all $T\/$ is said to be {\em repetitive}.%
\footnote{There is a parallel concept for tilings, where various terms are
used: ``almost periodic'' in \cite{Sol97}, ``tilings with local isomorphism''
in \cite{GS87,RW92}, and ``repetitive'' in \cite{Sen95}.  For
symbolic dynamical systems the term ``uniformly recurrent'' is used.}

(ii) A Delone set $X\/$ is {\em linearly repetitive} if there
is a constant $c$ such that $M_X(T) < c T$ for all $T > 0.$
}
\end{defi}

Repetitive Delone sets necessarily have finite local complexity,
since the definition implies that $N_X(T)$ is finite for all $T$.
The repetitivity function $M_X (T)$ of a Delone set provides a second
quantitative measure of its complexity, supplementing $N_X(T).$
Ideal crystals are repetitive and have $M_X(T)$ bounded.

Our second main result is that except for ideal crystals the slowest
possible growth rate of $M_X(T)$ is at least linear in $T$.

\begin{Theorem}~\label{th22} 
If a Delone set $X$ in $\RR^n$ has a single value of $T>0$ such that
\beql{Mcrystal}
M_X(T) < \frac {1}{3}\;T,
\eeq
then $X$ is an ideal crystal.
\end{Theorem}

The constant $\frac{1}{3}$ here is independent of the Delone constants
$r\/$ and $R\/$, which is made possible by the fact that the repetitivity
function is scale-invariant, and it is even independent of the dimension $n$.
It cannot be increased to be larger than $\frac{1}{2}$ (for any $n\/$)
as we show in \S\ref{examples}.

For $n\ge2$ we can show that a slightly larger value of the constant
in Theorem~\ref{th22}, while not necessarily guaranteeing that $X\/$
is an ideal crystal, at least gives a sufficient condition for $X\/$
to have a non-zero period.  Before stating this result we need to
define some constants.

\begin{defi}~\label{defs.5}
{\em
(i) For a Delone set $X\/$ we call the ratio $\kappa_X =R/r\/$ the
{\em packing-covering ratio\/} of $X$.  

(ii) The {\em Delone packing-covering constant} $\kappa(n)$, for dimension
$n$, is the infimum of $\kappa_X$ over all Delone sets $X\/$ in $\RR^n$.

(iii) The {\em lattice packing-covering constant} $\kappa_L(n)$, for
dimension $n$, is the infimum of $\kappa_\La$ over all lattices $\La$
in $\RR^n$.
}
\end{defi}

These constants were intoduced  by Ryshkov~\cite{R74}, who actually studied
the quantities $\frac {1}{2}\kappa(n)$ and $\frac {1}{2}\kappa_L(n)$
(corresponding to taking $2r\/$ instead of $r\/$ as the first Delone constant,
as mentioned after Definition~\ref{defs.1}).
If we now put
$$
c(n) :=\frac{\kappa(n)}{\kappa(n)+2}
$$
we have:

\begin{Theorem}~\label{th23} 
If a Delone set $X$ in $\RR^n$ has a single value of $T > 0$ such that
\beql{Mperiod}
M_X(T) < c(n)\,T,
\eeq
then $X$ has a non-zero period.
\end{Theorem}

Clearly $\kappa(1)=1$, so $c(1)=\frac{1}{3}$ (reflecting the
fact that in dimension~1 having a non-zero period is equivalent
to being an ideal crystal), and $\kappa(n) > 1$ for $n\ge2$.
In \S7 we present some results about  $\kappa(n)$ and $\kappa_L(n)$, 
including
the  result of 
Ryshkov~\cite{R74}
that $\kappa(n) \le 2$, which implies
that $\frac{1}{3} < c(n) \leq \frac {1}{2}$ for $n \geq 2$.
Thus Theorem~\ref{th23} gives at most a very slight dimension-dependent
improvement on Theorem~\ref{th22}.
In Theorem~\ref{th62} in \S6 we show that one cannot increase $c(n)$
in Theorem~\ref{th23} to any value exceeding $\frac{1}{2}\kappa(n)\le1$. 

%
%

\section{Bounded Patch-Counts and Ideal Crystals}\label{Nbdd}
\hsp
A Delone set $X\/$ is an ideal crystal if and only if
its patch-counting function $N_X(T)$ is bounded.
We have the following more precise result.

\begin{Theorem}~\label{th31}
If a Delone set $X$ in $\RR^n$ has $N_X(T)$ bounded with maximum value
$N$, then $X$ is a union of $N$ cosets of a lattice in $\RR^n$.
\end{Theorem}

\paragraph{Proof.}For each $T>0$ we can classify the points of $X\/$
according to their $T\/$-patches.  If {\bx} and {\by} have different
$T\/$-patches then they have different $T'$-patches for every $T'>T$,
so increasing $T\/$ refines the classification.  If $N_X(T)$ has
maximum value $N\/$ there is a $U\/$ with $N_X(U)=N\/$ and classifying
points of $X\/$ according to their $U$-patches gives a partition
$$
X=X_1\cup X_2\cup\cdots\cup X_N
$$
of $X\/$ such that points in the same class have identical
$T\/$-patches for all $T\/$ and points in different classes have
different $T\/$-patches for all $T\ge U$.  There is no loss of generality
in translating $X\/$ so that ${\bf 0}\in X_1$.  For any two points
$\bx$ and $\by$ we have $B({\bf 0};U)\subseteq B(\by;U+\|\by\|)$ and
$B(\bx-\by;U)\subseteq B(\bx;U+\|\by\|)$.
When $\bx$ and $\by$ are both in $X_1$ their $(U+\|\by\|)$-patches are
equal, so the translation by $\bx-\by$ that takes $B(\by;U+\|\by\|)$
to $B(\bx;U+\|\by\|)$ maps the points of $X\/$ in $B({\bf 0};U)$
one-to-one onto the points of $X$ in $B(\bx-\by;U)$.
Thus $\bx-\by\in X$ and $\sP_X(\bx-\by;U)\sim\sP_X({\bf 0};U)$.
Hence $\bx-\by\in X_1$.  This shows that $X_1$ is an additive
subgroup of $\RR^n$.

Next take any $\bx_1\in X_1$ and $\bx_2\in X_2$ (if $N>1$).
Then $\sP_X(\bx_1;U+\|\bx_2\|)\sim\sP_X({\bf 0};U+\|\bx_2\|)$,
and in a similar way we can use the translation by $\bx_2$ to show
that $\bx_1+\bx_2\in X$ and has the same $U\/$-patch as $\bx_2$.
Hence $\bx_1+\bx_2\in X_2$ and $X_2\supseteq X_1+\bx_2$.
Finally, for any other point $\by_2\in X_2$ we have
$\sP_X(\by_2;U+\|\bx_2\|)\sim\sP_X(\bx_2;U+\|\bx_2\|)$ and the
translation by $-\bx_2$ shows that $\by_2-\bx_2\in X_1$.
Hence $X_2-\bx_2\subseteq X_1$ so $X_2=X_1+\bx_2$ is a coset
of $X_1$ in $\RR^n$.
In the same way, each $X_i$ is a coset of $X_1$.

Since the Delone set $X\/$ is the union of finitely many translates
of $X_1$, $X_1$ itself is a Delone set and, being a subgroup of
$\RR^n$, is a full rank lattice.$~~~\Box$

\begin{Corollary}\label{cor31}
For Delone sets $X$,
\begin{eqnarray*}
\mbox{$N_X(T)=1$ for all $T$}&\Longleftrightarrow&
\mbox{$X\/$ is a translate of a lattice, and}\\
\mbox{$N_X(T)$ is bounded}&\Longleftrightarrow&
\mbox{$X$ is an ideal crystal.}
\end{eqnarray*}
\end{Corollary}
\paragraph{Proof.}
The implications in one direction follow from Theorem~\ref{th31}.
For the other direction, if $X=\La+F\/$ is an ideal crystal then the
orbits of $X\/$ under its group of translation symmetries are cosets
of a lattice containing the lattice $\La$ and the number of orbits
is at most the cardinal of the finite set $F$.  Clearly points in
the same orbit have identical $T$-patches for all $T\/$.$~~~\Box$

%
%

\section{Linear Patch-Counts and Ideal Crystals}\label{Nlinear}
\hsp
The object of this section is to show that if $X\/$ is not an
ideal crystal then there is a linear lower bound on the growth rate 
of $N_X(T)$.
The main step in doing this is the following lemma, which is analogous
to the result in the combinatorics of words that if the word-counting
function of an infinite sequence of symbols is the same for two
consecutive word lengths then it remains the same for all greater
word lengths (cf.\ our proof of Theorem~\ref{th41} below).

\begin{Lemma}\label{Lemma}
If $X$ is a Delone set with covering radius $R\/$ and there
are radii $U>0$ and $V>U+2R\/$ such that $N_X(V)=N_X(U)$ then
$N_X(T)=N_X(U)$ for all $T\ge U$, and $X\/$ is an ideal crystal.
\end{Lemma}

\paragraph{Proof.}If $\bx$ and $\by$ are points of $X\/$ with
$\sP_X(\bx;U)\not\sim \sP_X(\by;U)$ then clearly
$\sP_X(\bx;V)\not\sim\sP_X(\by;V)$.
So if $N_X(V)=N_X(U)$ then, for each $\bx\in X$, $\sP_X(\bx;V)$
is determined uniquely by $\sP_X(\bx;U)$.

We next show that $\sP_X(\bx;2V-U-2R)$ is also determined uniquely by
$\sP_X(\bx;U)$.  Take any $\bx_0\in X$, fixed for the moment, and
consider a point $\bc\in B(\bx_0;2V-U-2R)$.  We can find a point
$\bd\in B(\bx_0;V-U-R)$ with $\|\bd-\bc\|<V-R\/$ and, by the definition of
$R$, a point $\bx$ of $X\/$ in $B(\bd;R)$.  Then 
$B(\bx;U)\subset B(\bx_0;V)$
and $\bc\in B(\bx,V)$.  So $\sP_X(\bx_0;V)$ determines 
$\sP_X(\bx;U)$ which
in turn determines $\sP_X(\bx;V)$.  
Since $\bc$ was an arbitrary point in $B(\bx_0;2V-U-2R)$,
$$
B(\bx_0;2V-U-2R)\subset\bigcup_\bx B(\bx;V)
$$
where the union is over all points $\bx\in\sP_X(\bx_0;V-U)$, 
so the points
of $X\/$ in $B(\bx_0;2V-U-2R)$ are determined uniquely by 
$\sP_X(\bx_0;V)$.
Since $\bx_0$ was an arbitrary point of $X$,
$$
N_X(2V-U-2R)=N_X(V)=N_X(U).
$$
Iteration of this argument successively increases $V\/$ by $V-U-2R\/$
and shows that $N_X(T)=N_X(U)$ for all $T\ge U$.$~~~\Box$

Theorem~\ref{th21} is an almost immediate consequence of this lemma.

\paragraph{Proof of Theorem~\ref{th21}.}
If $X\/$ is not an ideal crystal then, by Corollary~\ref{cor31}, 
$N_X(T)$ is unbounded.
Certainly  $N_X(\eps_0)\ge 1$ for any $\eps_0>0$.  
Since $N_X(T)$ is unbounded
$$N_X(2R+\eps_0+\eps_1) > N_X(\eps_0) \mbox{ for any }\eps_1>0,$$
by Lemma~\ref{Lemma},
where $R\/$ is the covering radius of $X$.  Hence $N(2R+\eps_0+\eps_1)\ge2$.
Repeating this argument we find $N(4R+\eps_0+\eps_1+\eps_2)\ge3$, \dots,
$N(2mR+\eps_0+\eps_1+\cdots+\eps_m)\ge m+1$ for any integer $m\ge0$,
where $\eps_0,\ldots,\eps_m$ can be chosen arbitrarily small.
It follows that $N_X(T)\ge T/2R\/$ for every $T>0$.$~~~\Box$

One can prove an analogous result in multi-dimensional symbolic
dynamics, giving a condition for full periodicity of a symbol
pattern on the lattice $\ZZ^n$.  We consider symbolic words $S$
on the lattice $\ZZ^n$, drawn from a finite alphabet
$\sA$, i.e.\ $S \in \sA^{\ZZ^n}.$ 

\begin{Theorem}~\label{th41}
Let $S$ be a symbolic word on the lattice $\ZZ^n$ in which $A$
different symbols occur and let $N_S(m)$ count the number of
different symbol patterns in $S\/$ over all lattice cubes with
$m$ lattice points along each side.  If, for some integer $m \geq 1$, 
\beql{407}
N_S(m) < m + A - 1,
\eeq
then $S\/$ is fully periodic, with a full rank period lattice
$\La\subset \ZZ^n$.
\end{Theorem}

\noindent\paragraph{Proof.}Suppose that $N_S(m_0+1)=N_S(m_0)$ for some $m_0$.
Then whenever two cubes of side $m_0$ with centers at $\bx$ and $\by$ (where
$\bx$ and $\by$ are lattice points if $m_0$ is odd and half lattice points
if $m_0$ is even) have identical symbol patterns so do the cubes of side
$m_0+1$ with centers at $\bx+(\frac{1}{2},\frac{1}{2},\ldots,\frac{1}{2})$
and $\by+(\frac{1}{2},\frac{1}{2},\ldots,\frac{1}{2})$, since otherwise
the number of symbol patterns occurring in cubes of side $m_0+1$ would
be at least one greater than the number occurring in cubes of side $m_0$.
The same holds when $(\frac{1}{2},\frac{1}{2},\ldots,\frac{1}{2})$ is
replaced by $(\pm\frac{1}{2},\pm\frac{1}{2},\ldots,\pm\frac{1}{2})$
for any choice of signs.
Consequently the symbol pattern on any cube of side $m_0$ determines
uniquely the symbol pattern on the cube of side $m_0+2$ with the same center.
(This cube consists of the original cube surrounded by an extra shell of
lattice points.)
Hence $N_S(m_0+2)=N_S(m)$ and it follows by induction that the symbol pattern
on any cube of side $m_0$ determines the pattern on the whole of $\ZZ^n$.
Take any vector $\mbox{\boldmath$\la$}\in\ZZ^n$.
Since there are only finitely many symbols there must be two cubes of
side $m_0$ with the same symbol pattern such that one is a translate
of the other by $p\mbox{\boldmath$\la$}$, for some $p\in\ZZ$.
Since each cube determines the same global pattern,
$p\mbox{\boldmath$\la$}$ is a period of $S$.
Since it has periods in every lattice direction, $S\/$ is fully periodic.

It follows that if $S\/$ is not fully periodic then $N_S(m)$ is a strictly
increasing function of $m$, and hence
$$
N_S(m)\ge N_S(1)+m-1=A+m-1.
$$
$~~~\Box$ 

\noindent\paragraph{Remark.}The one-dimensional case of
Theorem~\ref{th41} is well-known and is due to Morse and
Hedlund \cite[Theorem~7.4]{MH38} in 1938.
The proof above is an extension of the one-dimensional proof.
The result is best possible in all cases: one-dimensional Sturmian
sequences $S\/$ have $A=2$ and $N_S(m)=m+1$ for all $m>1$, and examples
in higher dimensions are given by patterns that are Sturmian in one
coordinate direction and constant in all other coordinate directions.
There exist optimal examples for all other values of $A\/$ too.

%
%

\section{Linear Repetitivity and Ideal Crystals}\label{M}
\hsp
In \cite{LP99} we studied {\em linearly repetitive\/} sets,
which are Delone sets $X\/$ for which there is some constant
$c\/$ with $M_X(T)<cT\/$ for all $T$.
There are many examples of aperiodic linearly repetitive sets
associated to self-similar constructions and in \cite{LP99} such
sets were proposed as models for ``perfectly ordered quasicrystals.''
Theorem~\ref{th22} asserts that linear repetitivity, with a
sufficiently small constant, is sufficient to force a set $X\/$
to be an ideal crystal, and Theorem~\ref{th23} asserts that linear
repetitivity with a slightly larger constant, depending on the
dimension, forces at least one period.

\paragraph{Proof of Theorem~\ref{th22}.}
Without loss of generality we may suppose that ${\bf 0}\in X$.
Consider the patch $\sP := \sP_X({\bf 0};T) = B({\bf 0};T) \cap X$,
where $T\/$ satisfies (\ref{Mcrystal}), and let
$\Sig_\sP := \{\by\in X:\sP_X(\by;T)\sim\sP\}$ be the set of
centers of $T$-patches of $X\/$ translation-equivalent to $\sP$.

We first note that any vector
$\bp\in\Sig_\sP\cap B({\bf 0};\frac{2}{3}T)$ is a period of $X$.
This is because $\bx \in \sP$ and $\bp \in \Sig_\sP$ implies
$\bx+\bp\in\sP+\bp=\sP_X(\bp;T)\subset X$.
If $\bx\in\sP_X({\bf 0};\frac{1}{3}T)$, then  
$\bx+\bp\in B({\bf 0};T) \cap X = \sP$ and so, for all $\by\in\Sig_\sP$,
$$
\textstyle
\bx\in\sP_X(\by;\frac{1}{3}T)\quad\Rightarrow\quad
\bx+\bp\in\sP_X(\by;T)\subset X.
$$
Since $M_X(T)<\frac{1}{3}T$ the balls $B(\by;\frac{1}{3}T)$ with
$\by\in\Sig_\sP$ cover $\RR^n$ and therefore $\bp$ is a period of $X$.

It remains to show that $B({\bf 0};\frac{2}{3}T)$ contains $n\/$ linearly
independent vectors $\bp_1,\bp_2,\ldots,\bp_n$ of $\Sig_\sP$.  Let $B_1$
be any ball of radius $\frac{1}{3}T$ with {\bf 0} on its boundary.
Since $\frac{1}{3}T>M_X(T)$, $B_1$ contains a non-zero point $\bp_1$ of
$\Sig_\sP$.  When $\bp_1,\ldots,\bp_i$ have been chosen let $B_{i+1}$
be any ball of radius $\frac{1}{3}T$ that touches the linear subspace
$\langle\bp_1,\ldots,\bp_i\rangle$ at {\bf 0}.  Then there is a vector
$\bp_{i+1}$ of $\Sig_\sP$ in $B_{i+1}$, which is necessarily linearly
independent of $\bp_1,\ldots,\bp_i$.  Continuing in this way we find
$n\/$ linearly independent vectors $\bp_1,\ldots,\bp_n$ which all lie in
$B({\bf 0};\frac{2}{3}T)$, since each $B_i\subset B({\bf 0};\frac{2}{3}T)$.
$~~~\Box$

\paragraph{Proof of Theorem~\ref{th23}.}Let $\sP$ and $\Sig_\sP$ be as in
the proof of Theorem~\ref{th22}.  Since $M_X(T)$ is at least as large as the
covering radius of $\Sig_\sP$ the argument used in the proof of that theorem
shows that if there are two distinct points $\bx,\by\in\Sig_\sP$ with
\beql{pf23}
\|\bx-\by\|<T-M_X(T)
\eeq
then $\bx-\by$ is a period of $X$.  But the packing radius of
$\Sig_\sP$ is $\le M_X(T)/\kappa(n)$, so we can certainly
find $\bx,\by\in\Sig_\sP$ satisfying (\ref{pf23}) provided
$$
2M_X(T)/\kappa(n)<T-M_X(T),
$$
that is, provided (\ref{Mperiod}) holds.$~~~\Box$

%
%

\section{Non-crystals of Low Patch Complexity}\label{examples}
\hsp
Here we give examples that limit the extent to which the constants
in Theorems~\ref{th21}, \ref{th22} and \ref{th23} can be increased.
We first point out that, for trivial reasons, the constant $\frac{1}{2R}$
on the right of (\ref{Ncrystal}) cannot be improved and the constant
$\frac{1}{3}$ on the right of (\ref{Mcrystal}) cannot be increased
beyond $\frac{1}{2}$ for general Delone sets.
This is because for any Delone set $Y\/$ with Delone constants $r,R\/$ we
have $N_Y(T)=1=\frac{1}{2r}T\/$ and $M_Y(T)=R=\frac{R}{2r}T\/$ when $T=2r$.
When $n=1$ there are non-crystals with $\frac{R}{r}$ arbitrarily
close to 1, so $\frac{1}{2r}$ can be made as close as we like to
$\frac{1}{2R}$ and $\frac{R}{2r}$ as close as we like to $\frac{1}{2}$.
Higher dimensional examples can be found by noting that if
$X=Y\times\eta\ZZ^{n-1}$ then $X\/$ is a non-crystal when $Y\/$
is and has $N_X(T)=N_Y(T)$ and $M_X(T)$ very close to $M_Y(T)$,
for this $T$, when $\eta$ is small.

The growth rates of $N_X(T)$ and $M_X(T)$ as $T\to\infty$ are more
fundamental characteristics of $X\/$ than the sizes of these functions
at a particular value of $T$, however, so we give below
examples that combine near-minimal sizes of $N_X(T)$ and $M_X(T)$
for a particular $T\/$ with small linear growth rates as $T\to\infty$.
In particular, these examples are linearly repetitive, which
implies that they have a great deal of regularity
in the form of the existence of uniform patch frequencies and of
an approximately linear address map, as shown in  
\cite[Theorems 6.1,  7.1]{LP99}.
Most of our examples are not only Delone sets but Meyer sets,
where a {\em Meyer set\/} is a Delone set $X\/$ whose set of
differences $X-X\/$ is also a Delone set, see \cite{La96}, \cite{Moo97}.

We begin with one-dimensional examples and then extend them to
higher dimensions by the technique of taking
the direct product with a lattice of small mesh.

\begin{Theorem}~\label{th60}
For each  $\epsilon > 0$ there exists a one-dimensional aperiodic
Meyer set $Y=Y(\eps)$ which is linearly repetitive
and has the following properties.

(i) There is a $T^* > 0$ with
\beql{601}
N_Y(T^*) < \Bigl(\frac{1}{2} + \epsilon\Bigr) \frac {T^*}{R}\quad
\mbox{and}\quad M_Y(T^*) < \Bigl(\frac{1}{2} + \epsilon\Bigr) T^*.
\eeq

(ii) There is a $T_0>0$ such that, for all $T>T_0$,
\beql{602}
N_Y(T) < (1 + \epsilon) \frac {T}{R}\quad
\mbox{and}\quad M_Y(T) < (\tau +  1 + \epsilon) T,
\eeq
where $\tau = \frac{1}{2}(1 + \sqrt{5}) \approx  1.6180.$
\end{Theorem}

\noindent\paragraph{Proof.}
(i) As described above, to satisfy (i) it is sufficient to construct a set
$Y\/$ with $\kappa_Y = \frac{R}{r} < 1 + 2\epsilon$ and to take $T^*=2r$.
We shall take for $Y\/$ a set of the form 
$$
Y=Y(\alpha;\de) = \{m + \de\langle m\alpha\rangle : m \in \ZZ \},
$$
where $\alpha$ is irrational, $|\de| < \frac{1}{2}$, and
$\langle x\rangle=x-\lf x\rf$ is the fractional part of $x$.
This set is {\em almost linear}, in the sense of \cite{La99a},
so is a Meyer set by Theorem~5.1 of \cite{La99a}.
The irrationality of $\alpha$ ensures that $Y\/$ is aperiodic.
The packing radius $r_Y$ of $Y\/$ satisfies
$$
r_Y \geq \frac{1}{2} ( 1 - |\de|)
$$
while the covering radius $R_Y$ satisfies
$$
R_Y \leq \frac{1}{2} ( 1 + |\de|),
$$
so 
$$
\kappa_Y = \frac{R_Y}{r_Y} \leq  \frac{1+|\de|}{1-|\de|},
$$
which is $<1+2\eps$ when $|\de|<\frac{1}{2}\eps$.

(ii) To bound $N_Y(T)$ for large $T\/$ we note that the length of
the $m$th interval of $Y=Y(\alpha;\de)$ is
$$
1+\de(\langle(m+1)\alpha\rangle-\langle m\alpha\rangle)=
1+\de\alpha-\de(\lf(m+1)\alpha\rf-\lf m\alpha\rf),
$$
so the intervals of $Y\/$ are of two lengths and the arrangement
of these lengths corresponds to the arrangement of the terms of
the Beatty sequence (or Sturmian word)
\beql{Beatty}
\{\lf(m+1)\alpha\rf-\lf m\alpha\rf\}.
\eeq
For investigating the symbolic dynamics of $Y$, the two lengths
can be replaced by two abstract symbols, $A\/$ and $B\/$ say.
The study of such symbolic sequences goes back to
the eighteenth century but the emergence of a coherent body
of theory dates from two papers of Morse and Hedlund~\cite{MH38,MH40}
in 1938 and 1940.
Coven and Hedlund~\cite{CH73} showed that for $m\ge1$ there
are exactly $m+1$ distinct words of length $m\/$ in such a sequence
when $\alpha$ is irrational.  Since the number of intervals that are
contained in a given $T$-patch of $Y$, or overlap it sufficiently
for the overlap to determine their length, is at most%
\footnote{Perhaps more transparently, we could use instead the upper
bound $\lceil T+|\de|\rceil$ for the total number of intervals that
meet the $T$-patch, at the cost of slightly larger constant terms in
(\ref{NY}) and (\ref{MY}).}
$2\lf T+|\de|\rf$, this result shows that
\beql{NY}
N_Y(T) \le 2\lf T+|\de|\rf+1
\le (1+|\de|)\frac{T}{R}+1+2|\de|
\eeq
for all $T>0$.

The analogue in symbolic dynamics of the repetitivity function is
the recurrence function, $R(m)$, introduced by Morse and Hedlund.
As in \S1, we work in terms of the function $M_S(m)=R_S(m)-m$,
the maximum distance between the leading symbols of two successive
identical words of length $m\/$ in a symbolic sequence.
It is proved in \cite[p.2]{MH40} (see also \cite{AB98}, for example)
that if $\{q_k\}$ is the sequence of denominators of the continued
fraction convergents of $\alpha$ and $q_{k-1}\le m <q_k$ then 
$M(m)=q_k+ q_{k-1}$.
If the continued fraction has bounded partial quotients (as happens for
quadratic irrationals, for example) then $q_k+q_{k-1}<(b+2)q_{k-1}$,
where $b\/$ is an upper bound for the partial quotients, so the symbolic
sequence (\ref{Beatty}) is linearly repetitive and hence $Y\/$ is too.
To get the bound in (\ref{602}) we take $\alpha=\tau$, so that $b=1$
and $q_k=F_k$, the $k$th Fibonacci number.
Recalling that the number of identifiable intervals in any $T$-patch
is at most $2\lf T+|\de|\rf$, we then have
\beql{FibM}
2M_Y(T)\le F_k+F_{k-1}+|\de|
\eeq
when
\beql{range}
F_{k-1}\le2\lf T+|\de|\rf<F_k,
\eeq
where the last term on the right of (\ref{FibM}) is to take account of the
fact that the distance between two points of $Y\/$ differs from the number
of intervals separating them by at most $|\de|$.  Using the identity
$$
F_k=\tau F_{k-1}+\Bigl(\frac{-1}{\tau}\Bigr)^{k-1}
$$
and the lower bound in (\ref{range}) we obtain
\beql{MY}
M_Y(T)\le(\tau+1)T+\frac{1}{2}+\Bigl(\tau+\frac{3}{2}\Bigr)|\de|
\eeq
for all $T>0$.
If we choose $|\de|<\frac{1}{2}\eps$ and $T_0>2(\eps^{-1}+\eps)$
then (\ref{NY}) and (\ref{MY}) give (\ref{602}). $~~~\Box$ 

\paragraph{Remark.}
We note in passing that the set $Y(\tau;-1/\sqrt5)$ is a 
is a scaled version of the well known Fibonacci quasicrystal.
It can be shown to satisfy $M_Y(T)\le\tau^2T+\tau^2/2\sqrt5$,
where the constant term is half the length of the long interval.
In this case it can be shown that this  bound for $M_Y(T)$
is optimal in the sense that there exist arbitrarily large
$T\/$ for which it is achieved.

It is now straightforward to extend Theorem~\ref{th60} to the
$n$-dimensional case by taking direct products with lattices of small mesh.

\begin{Theorem}~\label{th61}
In $\RR^n$, for each $\epsilon > 0$ there exists a non-crystalline
Meyer set $X= X(\eps)$ which is linearly repetitive and has the
following properties.

(i) There is a $T^*>0$ with
\beql{603}
N_X(T^*) < \Bigl( \frac{1}{2} + \epsilon \Bigr) \frac {T^*}{R}\quad
\mbox{and}\quad M_X(T^*) < \Bigl( \frac{1}{2} + \epsilon  \Bigr) T^*.
\eeq

(ii) There is a $T_0>0$ such that, for all $T>T_0$,
\beql{604}
N_X(T) < (1 + \epsilon) \frac {T}{R}\quad
\mbox{and}\quad M_X(T) < (\tau+1+\epsilon) T.
\eeq
\end{Theorem}

\noindent\paragraph{Proof.}
We take $Y = Y(\eps)$ to be the one-dimensional
set constructed in Theorem~\ref{th60} and
$$
X : =Y\times\eta\ZZ^{n-1},
$$
where $\eta$ is a sufficiently small positive constant.
For any $\eta$, $X\/$ is an $n$-dimensional linearly repetitive
Meyer set which is a non-crystal, but has an $(n-1)$-dimensional
lattice of periods.
Because of these periods
\beql{Ns}
N_X(T)=N_Y(T)
\eeq
for all $T$.  We also have
\beql{Ms}
M_X(T)^2\le M_Y(T)^2+M_{\eta\ZZ^{n-1}}(T)^2=M_Y(T)^2+\frac{1}{4}\eta^2(n-1)
\eeq
for all $T$, and taking $T\le\min(2r_Y,\eta)$ in (\ref{Ms}) gives 
\beql{Rs}
R_X^2\le R_Y^2+\frac{1}{4}\eta^2(n-1),
\eeq
where $R_X$ is the covering radius of $X$.
Now (\ref{Ns}), (\ref{Ms}) and (\ref{Rs}), together with
Theorem~\ref{th60} (with $\frac{1}{2}\eps$ in place of $\eps$),
give (\ref{603}) and (\ref{604}) when $\eta$ is small enough,
with $T^*$ and $T_0$ as in Theorem~\ref{th60}.$~~~\Box$ 

\paragraph{Remark.} 
For the examples $X(\eps)$ the constants on the right of the
inequalities (\ref{603}) and (\ref{604}) simultaneously approach
$\frac{1}{2R}$, $\frac{1}{2}$, $1$ and $\tau+1$.
The first of these is optimal in all dimensions, by Theorem~\ref{th21},
and it is conceivable that the others are also individually optimal.
We note that if the second value, $\frac{1}{2}$, is optimal then,
since $c(n)\le\frac{1}{2}$ for all $n$, every Delone set $X\/$
satisfying the hypothesis of Theorem~\ref{th23} is not only periodic
but also an ideal crystal.

We end this section by deriving attainable values for the complexity
of aperiodic (rather than merely non-crystalline) Delone sets.
Note that, since $1 < \kappa(n) \le 2$ for $n \geq 2$, the constants in
(\ref{605}) below are larger than those in (i) of Theorem~\ref{th61}, as is
to be expected, but are also smaller than those in (ii) of Theorem~\ref{th61}.

\begin{Theorem}~\label{th62}
(i) In each dimension $n\geq1$ and for each $\eps>0$, there exists an
aperiodic Delone set $X=X(\eps)$ for which there is a $T>0$ such that
\beql{605}
N_X(T) < \Bigl( \frac{1}{2}\kappa(n)+ \epsilon \Bigr) \frac {T}{R}\quad
\mbox{and}\quad M_X(T) < \Bigl( \frac{1}{2}\kappa(n) + \epsilon  \Bigr) T.
\eeq

(ii) In each dimension $n \geq 1$ and for each $\epsilon > 0$, there
exists an  aperiodic repetitive Meyer set $Z=Z(\eps)$ having uniform
patch frequencies and for which there is a $T>0$ such that
\beql{606}
N_Z(T)<\Bigl(\frac{1}{2}\kappa_L(n)+\epsilon\Bigr)\frac {T}{R}\quad
\mbox{and}\quad M_Z(T)<\Bigl(\frac{1}{2}\kappa_L(n)+\epsilon\Bigr)T.
\eeq
\end{Theorem} 
\noindent\paragraph{Proof.}
(i) Take a Delone set $X\/$ with $\kappa_X<\kappa(n)+\eps$.
We may suppose that $X\/$ is aperiodic by moving one point of it by
a small enough amount not to disturb this inequality but in such a
way as to destroy any global periods.  Choose $T = 2r_X$, so that
$N_X (T) = 1 $ and $M_X (T) = R_X$.
Then (\ref{605}) holds since
$$
\Bigl(\frac{1}{2}\kappa(n)+\eps\Bigr)T>
\Bigl(\frac{R_X}{2r_X}+\frac{1}{2}\eps\Bigr)2r_X>R_X.
$$

(ii) Uniform patch frequencies are
defined as in \cite{LP99}. 
Take a lattice set $\La\in\RR^n$ with
$\kappa_\La<\kappa_L(n)+\eps$ and let
$$
Z=\{\mbox{\boldmath$\la$}+
\langle\phi(\mbox{\boldmath$\la$})\rangle\mbox{\boldmath$\de$}
:\mbox{\boldmath$\la$}\in\La\},
$$
where $\mbox{\boldmath$\de$}\in\RR^n$ and $\phi$ is a linear functional
on $\RR^n$.
By choosing {\boldmath$\de$} small enough (for a fixed $\phi$) we can
ensure that the packing and covering radii of $Z\/$ remain close enough to
the corresponding values for $\La$ for $\kappa_Z<\kappa_L(n)+\eps$ to hold.
Then, as in the proof of (i), (\ref{606}) is satisfied with $T=2r_Z$.
For the other properties of $Z\/$ we note that, up to an invertible
linear transformation, $Z\/$ is a cut-and-project set with lattice
$\La\times\ZZ\subset\RR^n\times\RR=\RR^{n+1}$, the
hyperplane $x_{n+1}=\phi(x_1,\ldots,x_n)$ as ``physical space''
and the line joining $(\mathbf{0},0)$ to
$(-\mbox{\boldmath$\de$},-1-\phi(\mbox{\boldmath$\de$}))$
as ``internal space'', the window being the half-open interval,
closed at the $(\mathbf{0},0)$ end, with these points as end points.
(The invertible linear transformation in question is projection
on the $\RR^n$ factor, its inverse being given by
$\mathbf{x}\mapsto(\mathbf{x},\phi(\mathbf{x})$.)
The facts that $Z\/$ is a repetitive Meyer set with uniform patch
frequencies that is aperiodic can now be derived from Propositions~2.21
and 2.23 of \cite{P} provided certain conditions are satisfied;
namely, (a) the physical space contains no non-zero vector of the
lattice $\La+\ZZ$, (b) the physical space is not contained in any
subspace generated by $n\/$ lattice vectors, and (c) the internal
space also contains no non-zero lattice vector.
Of these, (a) amounts to choosing $\phi$ so that
$1,\phi(\mbox{\boldmath$\la$}_1),\ldots,\phi(\mbox{\boldmath$\la$}_n)$
are linearly independent over $\QQ$, where
$\mbox{\boldmath$\la$}_1,\ldots,\mbox{\boldmath$\la$}_n$ is a basis of $\La$.
Then (b) holds too, since the physical space, having dimension $n$,
cannot lie in a subspace generated by $n\/$ lattice vectors unless
it contains those lattice vectors.
Finally, any lattice vector in internal space has the form
$(t\mbox{\boldmath$\de$},t(1+\phi(\mbox{\boldmath$\de$}))$ with
$t\mbox{\boldmath$\de$}\in\La$ and $t(1+\phi(\mbox{\boldmath$\de$}))\in\ZZ$,
so the internal space will contain no non-zero lattice vectors if
{\boldmath$\de$} is chosen to avoid the countably many vectors
$(m-\phi(\mbox{\boldmath$\la$}))^{-1}\mbox{\boldmath$\la$}$
with $\mbox{\boldmath$\la$}\in\La$, $m\in\ZZ$ and
$\phi(\mbox{\boldmath$\la$}),m\ne0$.$~~~\Box$

Note that while the set $Z\/$ in (ii) of Theorem~\ref{th62} is
repetitive, and hence also of finite local complexity, the set
$X\/$ in (i) may  have $N_X(T)$ or $M_X(T)$ infinite for large $T$.

To give the case $n=24$ as an example, using (\ref{n=24}) below;
by Theorem~\ref{th23} if $X\subset\RR^{24}$ satisfies
$$
M_X(T) <\frac{10\sqrt3-12}{13}T\approx 0.4093 T
$$
for some $T>0$ then $X\/$ has a non-zero period, but by
Theorem~\ref{th62}(ii) for any $c>1/\sqrt2\approx 0.7071$
there exist aperiodic sets $X\/$ in $\RR^{24}$ with uniform
patch frequencies that satisfy $M_X(T)< cT\/$ for some $T>0$.

%
%

\section{Delone Packing-Covering Constant}
\hsp
Recall that the {\em Delone packing-covering constant} $\kappa(n)$ in $n\/$
dimensions is the infimum of $\frac{R}{r}$ taken over all $n$-dimensional
Delone sets and the {\em lattice packing-covering constant} $\kappa_L(n)$
is the infimum of $\frac{R}{r}$ over all $n$-dimensional lattices.
In this section we prove various facts about $\kappa(n)$ and $\kappa_L(n)$.
We first show that these infimums are attained. 

\begin{Theorem}~\label{th71}
In each dimension $n$ there exists a Delone set $X$ whose Delone
constants $(r, R)$ satisfy $\frac{R}{r} =\kappa(n)$, and a lattice
$\La$ whose Delone constants satisfy $\frac{R}{r} = \kappa_L(n)$.
\end{Theorem}

\paragraph{Proof.} 
For each $m\in \NN$ let
$$
\sS_m=\{S\subset\RR^n:|\bx|\le m\mbox{ and }
|\bx-\by|\ge1\mbox{ for all }\bx,\by\in\sS\}.
$$
Then each $S\in\sS_m$ contains only finitely many points and
the Hausdorff metric
$$
d(S_1,S_2)=\max_{\bx\in S_1}\min_{\by\in S_2}|\bx-\by|+
\max_{\bx\in S_2}\min_{\by\in S_1}|\bx-\by|
$$
makes $\sS_m$ a compact metric space.  After suitable scalings we can
find a sequence $\{X_i\}$ of Delone sets in $\RR^n$ having $r\ge1$
for all $i\/$ and $R\/$ tending to $\kappa(n)$ as $i\to\infty$.  The sets
$X_i\cap B(\mbox{\boldmath$0$};1)$ belong to the compact space $\sS_1$,
so we can find a convergent subsequence
$$
\{S_{1i}\}=\{X_{1i}\cap B(\mbox{\boldmath$0$};1)\}
$$
which we can choose so that $d(S_{1i},S_{11})<1$ for all $i$.
Similarly the sets $X_{1i}\cap B(\mbox{\boldmath$0$};2)$ belong
to $\sS_2$, and we can find a subsequence
$$
\{S_{2i}\}=\{X_{2i}\cap B(\mbox{\boldmath$0$};2)\}
$$
with $d(S_{2i},S_{21})<1/2$ for all $i$.  Continuing in this way, for
each $j\in\NN^+$ we can find a subsequence $\{X_{ji}\}_{i=1}^\infty$
of $\{X_{j-1,i}\}_{i=1}^\infty$ such that
$$
d(X_{ji}\cap B(\mbox{\boldmath$0$};j),X_{j1}\cap B(\mbox{\boldmath$0$};j))
<\frac{1}{j}\quad\mbox{for all $i\in\NN^+$.}
$$
The sequence of Delone sets $\{X_{j1}\}_{j=1}^\infty$ now converges
to a unique set $X\subset\RR^n$ that has $r\ge1$ and $R\le\kappa(n)$,
and by the definition of $\kappa(n)$ there is equality in both places.

The proof for the lattice case is similar, noting that a limit of
lattices that is a Delone set is a lattice.
$~~~\Box$

We next give upper and lower bounds for $\kappa(n)$, which are
due to Ryshkov~\cite{R74}.

\begin{Proposition}~\label{prop71}
In each dimension $n\/$ the Delone packing-covering constant
$\kappa(n)$ satisfies
\beql{702}
\sqrt{\frac{2n}{n+1}}\le\kappa(n) \leq 2.
\eeq
\end{Proposition}

\noindent\paragraph{Proof.}
The lower bound for $\kappa(n)$ follows by combining Rogers's upper bound
for the packing density \cite[Theorem~7.1]{Z99} and Coxeter, Few and
Rogers's lower bound for the covering thickness \cite[Theorem~3.4]{Z99} of
equal balls in $\RR^n$.  The former asserts that for every $\eps>0$ the
number of balls of radius $r\/$ that can be packed into a ball of radius
$l\/$ is $<(1+\eps)\sig_n(l/r)^n$ when $l\/$ is large enough, where $\sig_n$
is the proportion of an $n$-dimensional regular simplex that is covered
by balls centered at its vertices and reaching to the mid-points of its
sides.  The latter asserts that the number of balls of radius $R\/$ it
takes to cover a ball of radius $l\/$ is $ > (1-\eps)\tau_n(l/R)^n$ when
$l\/$ is large enough, where $\tau_n$ is the proportion of the simplex
(counted with multiplicity) that is covered by balls centered at its
vertices and reaching to its centroid (these latter balls just cover
the simplex).  It follows from these bounds that the packing radius
$r\/$ and covering radius $R\/$ of any set $X\subset\RR^n$ satisfy
$$
\frac{R}{r}\ge\Biggl(\frac{\tau_n}{\sig_n}\Biggr)^{1/n}=
\sqrt{\frac{2n}{n+1}}\;,
$$
the expression on the right being the ratio of the circumradius to half
the edge length of a regular $n$-simplex.  

To establish the upper bound for $\kappa(n)$, take a packing or $\RR^n$
with unit spheres, to which no sphere can be added without overlap.  The
sphere centers then form a Delone set $X$ with $r \geq 1$.  We claim that
$R\leq2$. If not, and $R > 2$, then there would exist some point at a
distance at least $2$ from all sphere centers, and a unit sphere placed
at this point would overlap no sphere in the packing, a contradiction.
$~~~\Box$

The lower bound in Proposition~\ref{prop71} is sufficient to establish
that $\kappa(1)=\kappa_L(1) = 1$ and $\kappa(2)=\kappa_L(2)=2/\sqrt3$, 
with the optimal configurations in dimensions 1 and 2 being $\ZZ$ and
the hexagonal lattice $A_2$.
This gives $c_1=\frac{1}{3}$ and $c_2=\frac{1}{2}(\sqrt3-1)\approx0.3660$.
In three dimensions  we have the bounds
$$
\sqrt{\frac{3}{2}}\approx1.2247\le\kappa(3)\le \kappa_L(3)=
\sqrt{\frac{5}{3}}\approx1.2910,
$$
the upper bound coming from the body-centered cubic lattice $A_3^*$.
In twenty-four dimensions  we have the remarkably close bounds
\beql{n=24}
\frac{4\sqrt3}{5}\approx1.3856\le\kappa(24)\le\kappa_L(24)
\le\sqrt2\approx1.4142,
\eeq
with the upper bound coming from the Leech lattice.  (The values of
the packing density and covering thickness of the Leech lattice are
given in \cite{CS99}.)  Note that $\kappa(n)$ and $\kappa_L(n)$ are
probably not monotonically increasing in $n$. 

The lower bound in Proposition~\ref{prop71} also yields the asymptotic bound
$$
\liminf_{n \to \infty} \kappa(n) \geq \sqrt{2}.
$$
This asymptotic bound can be further improved by using the current best
upper bound for sphere packing density (the Kabatiansky-Levenshtein bound)
as $n \to \infty$, instead of the Rogers bound, which yields
$$
\liminf_{n \to \infty} \kappa(n) \geq 2^{0.5990} \approx 1.5146,
$$
cf.\ \cite[Chapters~1 and 9]{CS99}, \cite[Theorem~8.2]{Z99}.

The upper bound of Proposition~\ref{prop71} shows that the constant 
$c(n) = \frac{\kappa(n)}{\kappa(n) + 2} \leq \frac{1}{2}$, hence 
Theorem~\ref{th23} gives at most a very slight dimension-dependent
improvement on Theorem~\ref{th22}.

The lattice packing-covering constant $\kappa_L(n)$ satisfies the lower
bound in (\ref{702}) {\em a fortiori}.  In \cite{R74} Ryshkov
gives the upper bound
$$
\kappa_L(n)\le\sqrt\frac{n+2}{3}.
$$
(Recall that $\kappa_L(n)$ is, by definition, twice Ryshkov's constant.)
This is known to hold with equality in dimensions 1, 2 and 3.
There seems to be very little known  
about upper bounds for
$\kappa_L(n)$, and we wish to raise the following question.%
\footnote{Gruber and Lekkerkerker~\cite[p.618]{GL87} mistakenly
assert that Ryshkov~\cite{R74} proved the upper bound $\kappa_L(n)<2$;
in fact he proved only that $\kappa(n) \le 2$. Thus the question
raised here appears to be open.}

\noindent{\em Question. Is there a constant $C$ such that
$\kappa_L(n) \leq C$ for all $n$, or is $\kappa_L(n)$ unbounded
as $n \to \infty$?}

\noindent One reason for interest in the question is that in any
dimension $n$ for which $\kappa_L(n) > 2$ (should such dimensions
exist) the densest sphere packing with equal spheres cannot be a
lattice packing.

\noindent\paragraph{Acknowledgment.} We are indebted to
V. Berth\'{e} for comments and references 
concerning Theorem~\ref{th60}.

{\tt
\begin{tabular}{ll}
email: & jcl@research.att.com \\
  & pleasants\_p@usp.ac.fj
\end{tabular}
 }

\end{document}